\documentclass[preprint,12pt,square,numbers,sort&compress]{elsarticle}
\usepackage[normalem]{ulem}
\usepackage{graphicx,color} 
\usepackage{amsmath}  
\usepackage{amssymb}
\usepackage{epstopdf}
\usepackage{float}
\usepackage{subcaption}

\def\be{\begin{eqnarray}}
\def\ee{\end{eqnarray}}

\def\ds{\displaystyle}

\def\ni{\noindent}

\begin{document}

\title{Nonlinear second order inhomogeneous differential equations in one dimension}
%

\author{Yajnavalkya Bhattacharya \\ 
Jurij W. Darewych \\ York University, Toronto, Canada \\}

\begin{abstract}
  We study inhomogeneous nonlinear second-order differential equations in one dimension. The inhomogeneities can be point sources or continuous source distributions. We consider second order differential equations  of type $\phi''(x) + V(\phi(x)) = Q \, \delta(x) $, where $V(\phi)$ is a continuous, differentiable, analytic function and $Q \,\delta (x)$ is a point source. In particular we study cubic functions of the form $V(\phi(x)) = A\,\phi(x) + B\,\phi^3(x)$. We show that Green functions can be determined for modifications of such cubic equations, and that such Green's functions can be used to determine the solutions for cases where the point source is replaced by a continuous source distribution.      
\end{abstract}

\maketitle

\section{Introduction}

There are not many nonlinear differential equations in one dimension for which exact analytic solutions can be obtained in terms of elementary functions. This is especially so for nonlinear equations with an inhomogeneous source term.  This is noted in a textbook by Olver \cite{0.1}: \vskip .2cm \ni  
{\sl ``....the Superposition Principle for inhomogeneous linear equations allows one to combine the responses of the system to different external forcing functions...The two general Superposition Principles furnish us with powerful tools for solving linear partial differential equations, which we shall repeatedly exploit throughout this text. In contrast, nonlinear partial differential equations are much tougher, and, typically, knowledge of several solutions is of scant help in constructing others. Indeed, finding even one solution to a nonlinear partial differential equation can be quite a challenge.''}

Exact solutions to one-dimensional homogeneous nonlinear oscillators, such as the ``restricted Duffing'' equation
\be \label{mickens1}
\frac{d^2y}{dt^2 } + y + \epsilon y^3 = 0 \;\;\; \text{where}\;\;\;  \epsilon >0,\;\;\;y(0)=A,\;\;\;\frac{dy(0)}{dt}=0.
\ee
have been worked out in terms of Jacobi elliptic functions  by Mickens \cite{0.2a}, \cite{0.2b}.  

Approximate Green's functions for second order inhomogeneous nonlinear equations have been worked out    
by Frasca and Khurshudyan \cite{0.3}, and by others. For the nonrelativistic quartic oscillator, Anderson has derived the quantum mechanical Green's function $G(z_2, t_2; z_1, t_1)$ \cite{0.4}, based on an invertible linearization map \cite{0.5}.  

An analytic solution to the nonlinear (second-order) one-dimensional, cubic equation
\be \label{1} 
\left(-\frac{d^2}{d x^2 } - \mu^2\right) \, \phi(x) +  \lambda \, \phi^3(x) = 0 , 
\ee     
where $\mu$ and $\lambda$ are constants, is given by \cite{1}
\be \label{2}
\phi(x) = \frac{\mu}{ \sqrt\lambda} \tanh \left(\frac{x \mu}{\sqrt 2 } \right),            
\ee
as can be simply verified by substituting the function (\ref{2}) into eq. (\ref{1}).
\section{Inhomogeneous equations}
It turns out that the inhomogeneous generalisation of (\ref{1}), that is,  the nonlinear cubic eq. (\ref{1}) but with a point source 
on the right-hand side, namely  
\vskip .3 cm
\par \ni
\be \label{3} 
\left(-\frac{d^2}{d x^2 } - \mu^2\right) \, F(x) +  \lambda \, F^3(x) =  \delta (x), {\small {\scriptsize {\small }}}
\ee 
can be solved analytically.
This nonlinear, inhomogeneous eq. (\ref{3}) has the analytic solution:
\be \label{4}
F(x) = \frac{\mu}{ \sqrt\lambda} \tanh \left(\frac{|x| \mu}{\sqrt 2 } \right),
\ee 
This can be verified by substituting eq.(\ref{4}) into eq.(\ref{3}). Note that $|x|$ stands for $\sqrt{x^2}$, where, in one dimension, $-\infty<x<+\infty$.
Note also that the solution $F(x)$ of equation (\ref{4}) is just  $F(x) = \phi(|x|)$, where $\phi$ is given in eq. (\ref{2}).
\par \ni \vskip .2 cm
We plot the solutions (\ref{4}) in units of $\mu$ for a few values of $\lambda$ in figure \ref{fig:tanhmod}:
\begin{figure}[H]
\centering
\includegraphics[totalheight=9.1cm]{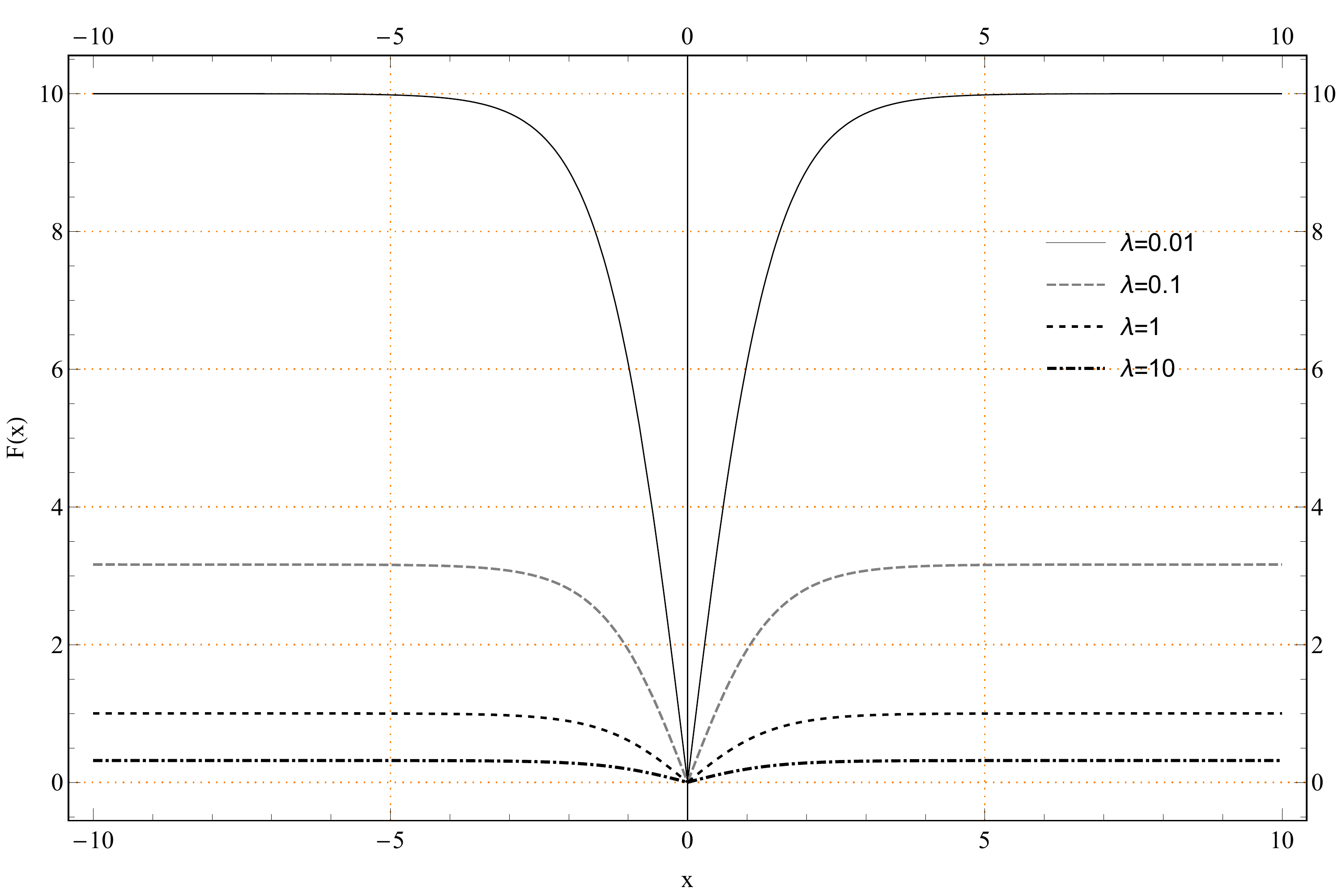}
\caption{Plot of $\ds F(x) = \frac{\mu}{ \sqrt\lambda} \tanh \left(\frac{|x| \mu}{\sqrt 2 } \right)$ in units of $\mu$=1, for various $\lambda$}
\label{fig:tanhmod}
\end{figure}
\vskip .5 cm
\par \ni
\vskip .3 cm
\par
\begin{figure}[H]
\centering
\includegraphics[totalheight=9.2cm]{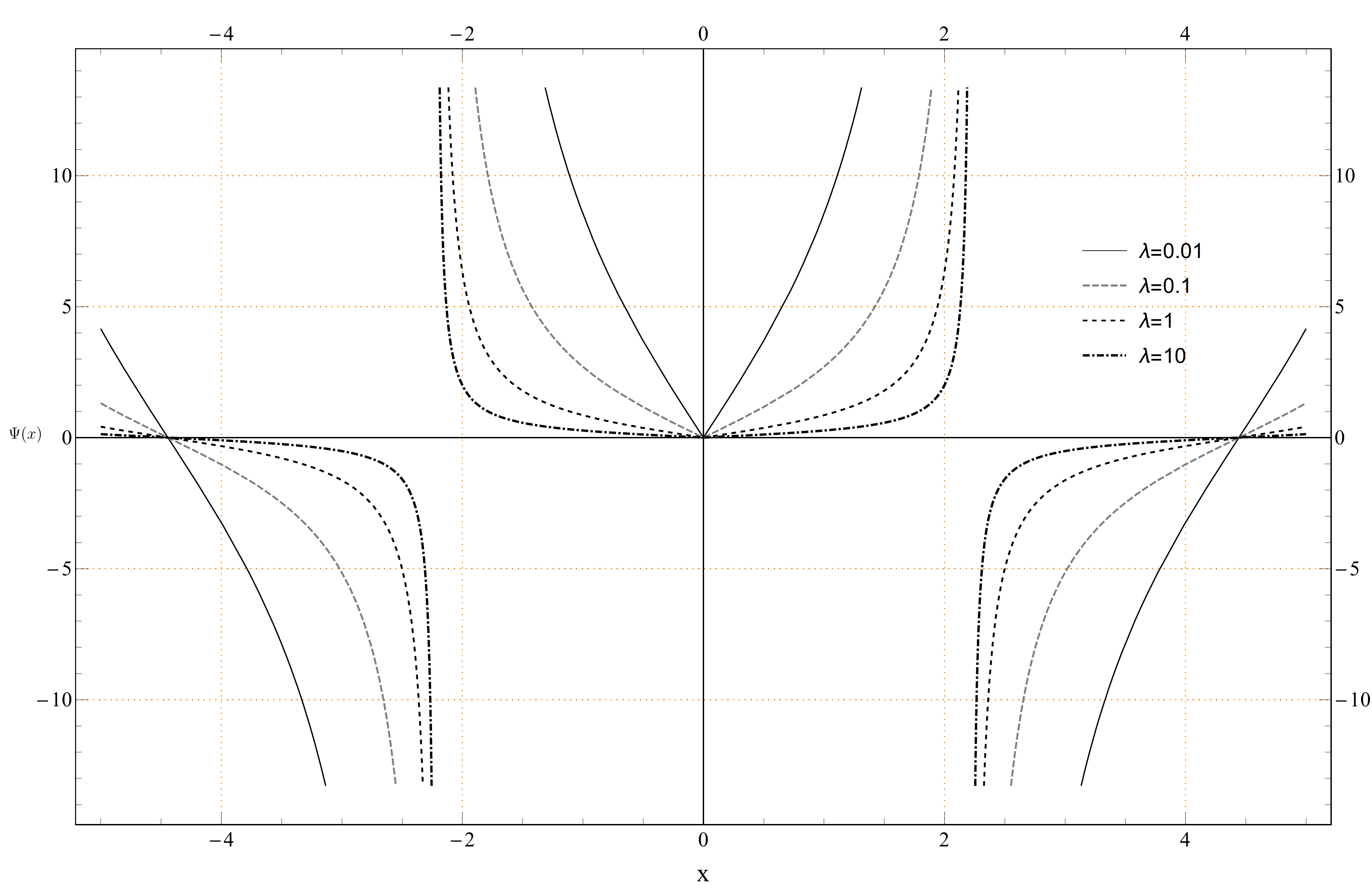}
\caption{Plot of $\ds \Psi(x) = \frac{\mu}{ \sqrt\lambda} \tan \left(\frac{|x|~m}{\sqrt 2 } \right)$ in units of $m$=1, for various $\lambda$}
\label{fig:tanmod}
\end{figure}

 \ni
If $\mu$ in equation (\ref{3}) is replaced by $i\,m$, where $m$ is also real like $\mu$, then the modified equation,  
\be \label{7} 
\left(-\frac{d^2}{d x^2 } + m^2\right) \, \psi(x) +  \lambda \, \psi^3(x) = 0 , 
\ee  
has the solution 
\be \label{8}
\psi(x) = \frac{m}{ \sqrt\lambda} \tan \left(\frac{x\;m}{\sqrt 2 } \right)                  
\ee
\ni
Thus, the inhomogeneous version of eq. (\ref{7}), 
\be \label{9} 
\left(-\frac{d^2}{d x^2 } + m^2\right) \, \Psi(x) +  \lambda \, \Psi^3(x) =  \delta (x), 
\ee 
has the solution 
\be \label{10}
\Psi(x) = \frac{m}{ \sqrt\lambda} \tan \left(\frac{|x| m}{\sqrt2} \right),                    
\ee
which is plotted in figure \ref{fig:tanmod}.
 \par \vskip .3 cm
\section{Generalisations}
We can generalise the above results for second-order equations to any equation of the form
\be \label{11}
- \frac{d^2 \phi(x)}{d x^2} + V\big(\phi(x)\big) = 0
\ee
where $V(\phi)$ is a continuous differentiable analytic function. 
 If the solution of eq. (\ref{11})  is $\phi(x)$, then the solution of the corresponding inhomogeneous equation with a point source, namely 
\be \label{12}
- \frac{d^2 \Phi(x)}{d x^2} + V\big(\Phi(x)\big) = \delta(x)
\ee
will be $\Phi(x) = \phi(|x|)$.
The second-order nonlinear equations (\ref{3}) and (\ref{9}),  as well as  their solutions (\ref{4}) and (\ref{10}), are particular examples of the general results (\ref{11}) and (\ref{12}). 
\par
A particularly simple form of (\ref{12}) is the linear  one-dimensional neutron diffusion equation with a plane source: 
\be \label {lininhomeq}
- \frac{d^2 \phi(x)}{dx^2} + k^2 \phi (x) =  \delta(x).
\ee
Eq. (\ref{lininhomeq}) corresponds to equation (\ref{12}) with $\ds V(x) = k^2 \, \phi (x)$. It has the solution 
\be \label {16}
\phi(x) = \frac{1}{2k} \, e^{-|kx|},
\ee
as can be verified by substituting eq. (\ref{16}) into eq. (\ref{lininhomeq}).
\par
Regardless of whether $\phi(x)$ of eq. (\ref{11}) is analytic, or can only be worked out numerically, it will nevertheless be true that the solution of (\ref{12}) will be $\phi(|x|)$.

 \section{Green's functions for the nonlinear equations and continuous source distributions}
 
 Green's function for the linear equation (\ref{lininhomeq})
 can be obtained by replacing the single point source $\delta(x)$ with a source of the form  $\delta (x-x_1)$.  Thus, the linear equation (\ref{lininhomeq})  becomes the equation 
 \be \label {17}
- \frac{d^2 \phi(x)}{dx^2} + k^2 \phi (x) = \delta(x-x_1),
\ee 
which has the solution  
 \be \label {18}
  \phi (x, x_1;k) = \frac{1}{2 k} e^{-k |x - x_1|},
 \ee
except for $x = x_1$ where it is singular. The result (\ref{18}) is recognized as the Green funtion for the Modified Helmholtz Equation in one dimension.\cite{5}
 \vskip .2 cm  \par
The analogous equation for the nonlinear case, cf. eq. (\ref{3}), is
 \be \label {19}
- \frac{d^2 G(x)}{dx^2} - \mu^2 G(x)  + \lambda \, G^3(x) = \delta(x-x_1).
\ee 
The solution, i.e. the Green function for the nonlinear eq. (\ref{19}), is 
\be \label{20}
G(x,x_1; \mu, \lambda) = \frac{\mu}{ \sqrt\lambda} \tanh \left(\frac{|x-x_1| \mu}{\sqrt 2 } \right),
\ee 
except at $x = x_1$ where $G$ is singular.

Analogously to (\ref{19}), eq. (\ref{9}) becomes 
\be  \label{21} 
\left(-\frac{d^2}{d x^2 } + m^2\right) \, \Psi(x) +  \lambda \, \Psi^3(x) =  \delta (x-x_1), 
\ee 
which has the solution, i.e. Green function,  
\be \label{Greenfunctionpsi}
\Psi(x,x_1; m, \lambda) = \frac{m}{ \sqrt\lambda} \tan \left(\frac{|x-x_1| m}{\sqrt2} \right).                    
\ee
\par \vskip .2cm
Once the Green functions are known, it is possible to calculate the ``potential" $V(x)$ due to a continuous source (``charge") distribution $R(x_1)$, by evaluating the integral:
\be \label{23} 
V(x) = \int_{-\infty}^{\infty} dx_1 \, R(x_1) \, G(x, x_1). 
\ee  
For the relatively simple linear case, cf.\,eqs (\ref{17}) and (\ref{18}), $V(x)$ can be evaluated analytically for various $R(x_1)$. For a Gaussian source distribution  $\ds R(x_1) = e^{- x_1^2} $, the potential $V_{\text lin}(x)$ 
\be \label{22} 
  V_{\text lin}(x) = \ds \int_{-\infty}^{\infty} dx_1 \, e^{- x_1^2} \, \frac{1}{2k} e^{-k |x - x_1|} ,
  \ee
evaluates to  
  \be \label{23a} 
V_{\text lin}(x)   = \ds \frac{\sqrt\pi}{4k} \,
e^{\left(-kx+\frac{k^2}{4}\right)}
\left[e^{2kx}Er\left(x+\frac{k}{2}\right) - e^{2kx} + Er\left(-x+\frac{k}{2}\right)-1 \right], 
\ee  
where  \:
$ \ds Er(u) = \frac{2}{\sqrt\pi} \int_{0}^{u} e^{- t^2} dt $ \;
is the Error function.
\vskip .2 cm \ni

\begin{figure}[H]
\centering
\includegraphics[totalheight=9cm]{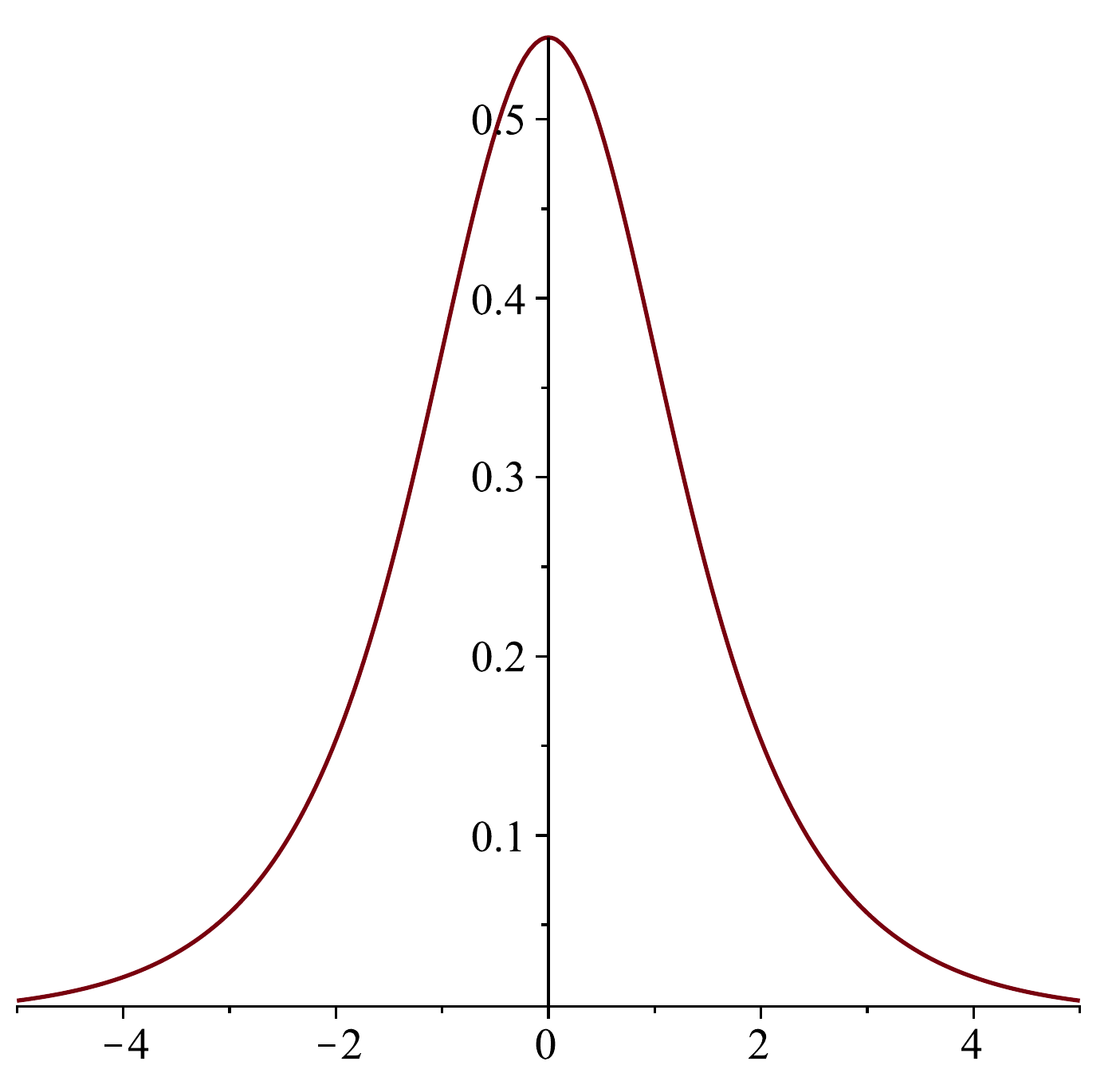}
\caption{Potential $\ds V_{\text lin}(x)$ = $\ds \int_{-\infty}^{\infty} dx_1 \, e^{- x_1^2} \, \frac{1}{2k} e^{-k |x - x_1|}$}.
\label{fig:Vxmaple}
\end{figure}
However, we are primarily concerned with the more interesting nonlinear cases (cf. eqs. (\ref{19}) and (\ref{21})). Unfortunately, the evaluation of $V(x)$, eq. (\ref{23}), for a given distribution $R(x_1)$, with either Green function (\ref{20}) or (\ref{Greenfunctionpsi}), generally must be done by numerical quadrature. 
\par
Analytic expressions for the pontential can be evaluated in a few cases; for example, a step-function source distribution, such as
\be \label{R1step}
R_1(x_1) = 1 \;\;\forall\;\;\{-4 \le x_1 \le 8\}, \;\;R_1 = 0\;\;\forall \{x_1>8,\;x_1<-4\}, 
\ee 
for which the potential is given by 
\be \label{vxstepfuncint}
\ds V_1(x) = \int_{-\infty}^{\infty} dx_1 \, R_1(x_1) \, G(x, x_1,1,1) 
\ee
The integral in eq. (\ref{vxstepfuncint}), which can be evaluated analytically by Maple\textsuperscript{\textregistered} and Mathematica\textsuperscript{\textregistered}, is plotted in Figure \ref{fig:stepfunction4to8}. The result is listed in Appendix I.   
\begin{figure}[H]
	\centering
	\includegraphics[totalheight=8cm]{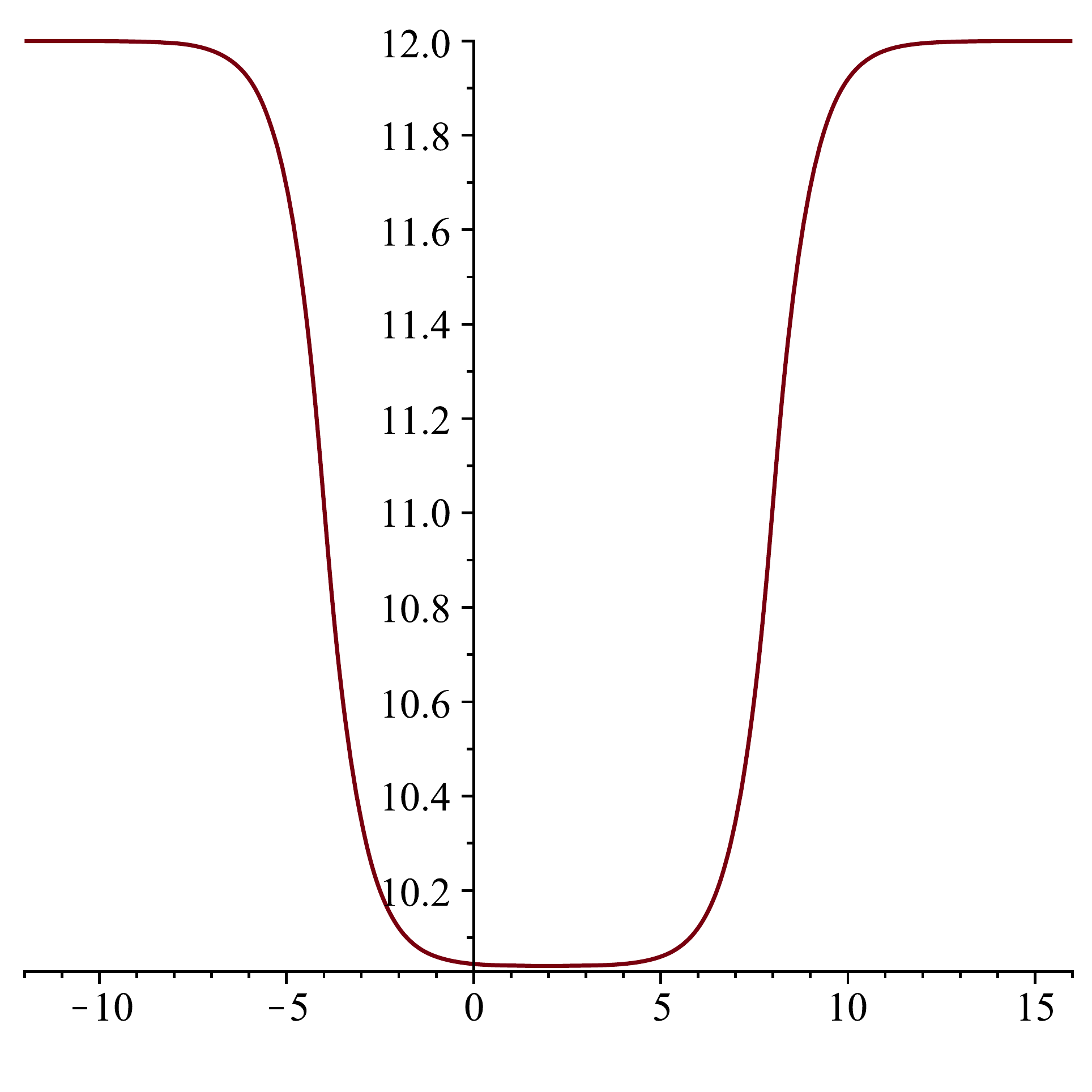}
	\caption{Plot of $\ds V_1(x)$ of eq. (\ref{vxstepfuncint}), for the ``step function'' distribution $R_1(x)$ of (\ref{R1step}).}
	\label{fig:stepfunction4to8}
\end{figure}

\begin{figure}[H]
	\centering
	\includegraphics[totalheight=7cm]{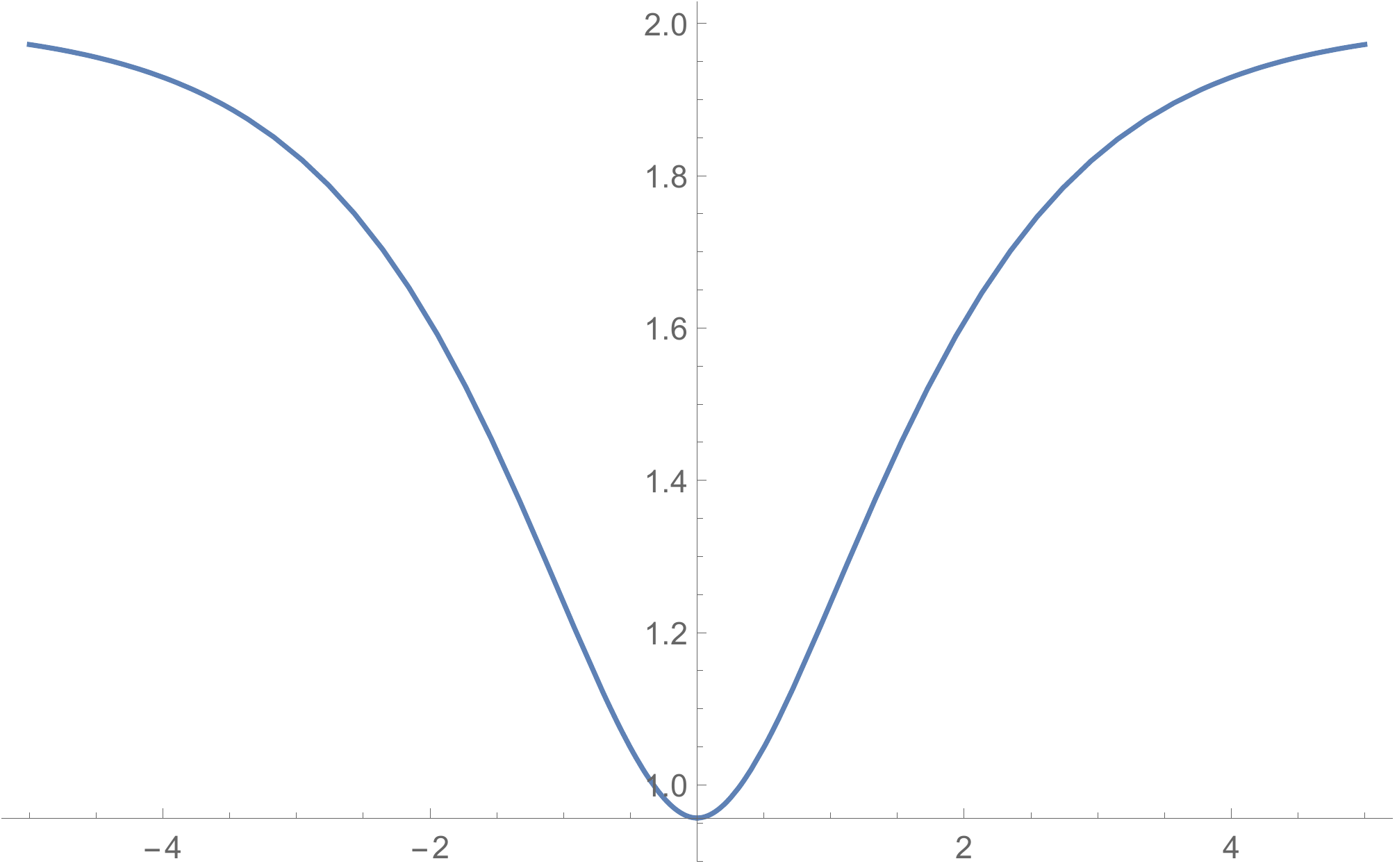}
	\caption{Plot of $\ds V_e(x)= \int_{-\infty}^{\infty} dx_1 \, R_e(x_1) \; G(x, x_1,1,1)\;$, for the ``exponential'' distribution $\ds R_e(x) = e^{-|x|}$. The plot of the potential $V_g(x)$ for the Gaussian distribution $e^{-x1^2}$ turns out to be  very similar in shape, and, thus, it is not plotted here. }
	\label{fig:exponentialmag}
\end{figure}

Two other distributions for which analytic solutions of eq. (\ref{23}) can be obtained are the ``exponential" distribution  $\ds R_e(x) = e^{-|x|}$, and the Gaussian distribution $R_g(x)=e^{-x1^2}$. The resulting potentials $\ds V_e(x) = \int_{-\infty}^{\infty} dx_1 \, R_e(x_1) \; G(x, x_1,1,1)$, and 
$\ds V_g(x) = \int_{-\infty}^{\infty} dx_1 \, R_g(x_1) \; G(x, x_1,1,1)$, though analytic, are cumbersome expressions. The solution for $V_e(x)$ is listed in Appendix I. Numerical integrations of the potentials yield identical results to their analytic counterparts, as should be expected.   
$V_{e}(x)$ is plotted in figure \ref{fig:exponentialmag}: 
\vskip .3 cm

\vskip .3 cm
\par 
Next, we exhibit examples of numerically calculated potentials $V(x)$ with the Green function of eq. (\ref{20}), for a few choices of source distributions, namely 

\be \label{gaussiansourcedist}
R_2(x_1) =  \ds e^{\ds-x_1^2/a^2}   
\ee 
\be \label{bellsourcedist}
R_3(x_1) = \frac{1} {\left[(x_1+a)^2+b^2)^2 \right]}
\ee
where $a$ and $b$ are arbitrary, real constants. The potentials are
\be \label{potential}
 V_j(x) = \int_{-\infty}^{\infty} dx_1\,R_j(x_1) \, G(x,x_1,1,1),
\ee
where $j = 2, 3$.
 
\begin{figure}[H]
 	\centering
 	\begin{subfigure}{0.48\textwidth}
 		\includegraphics[width=\textwidth]{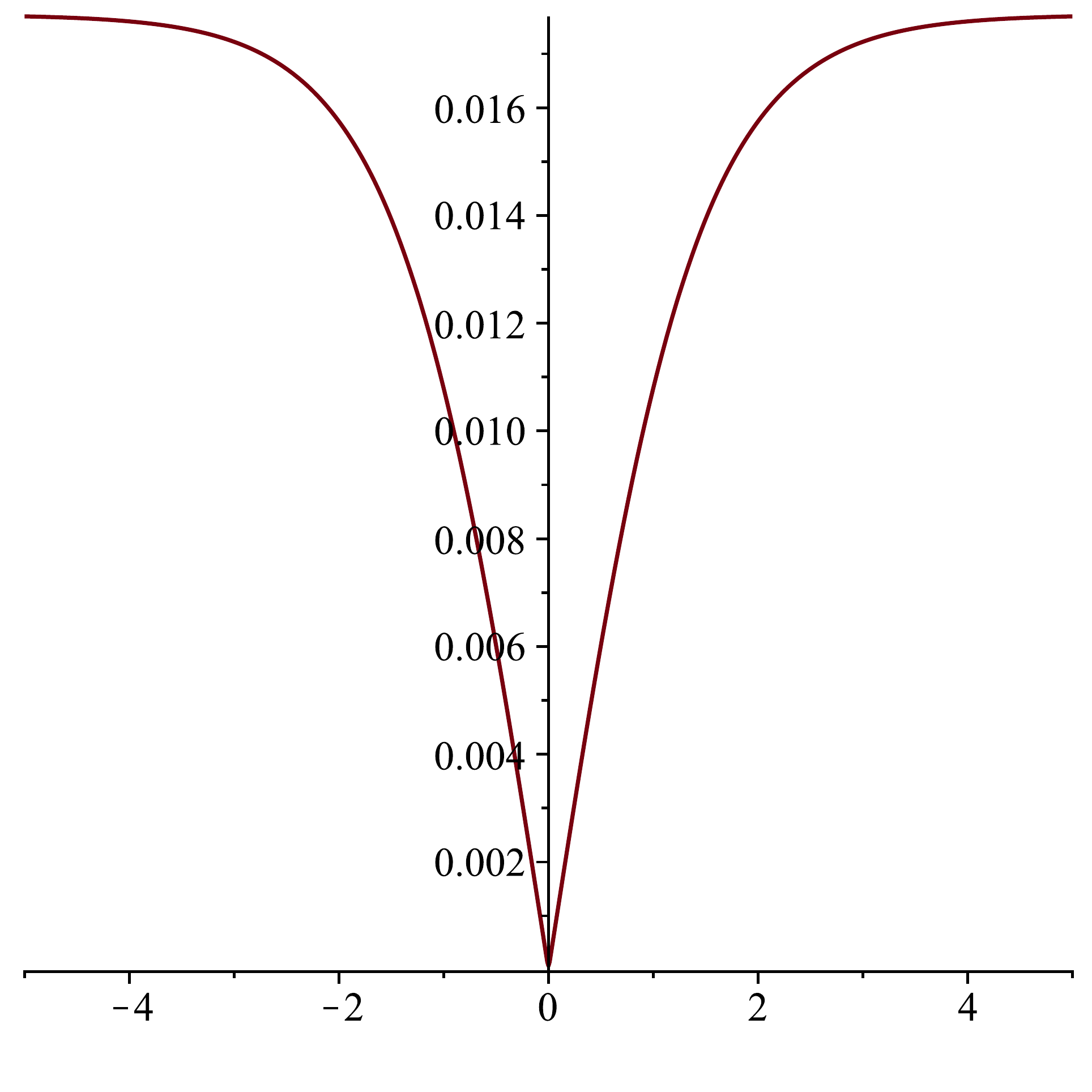}
 		\subcaption{}
 	\end{subfigure}
 	\begin{subfigure}{0.48\textwidth} 
 		\includegraphics[width=\textwidth]{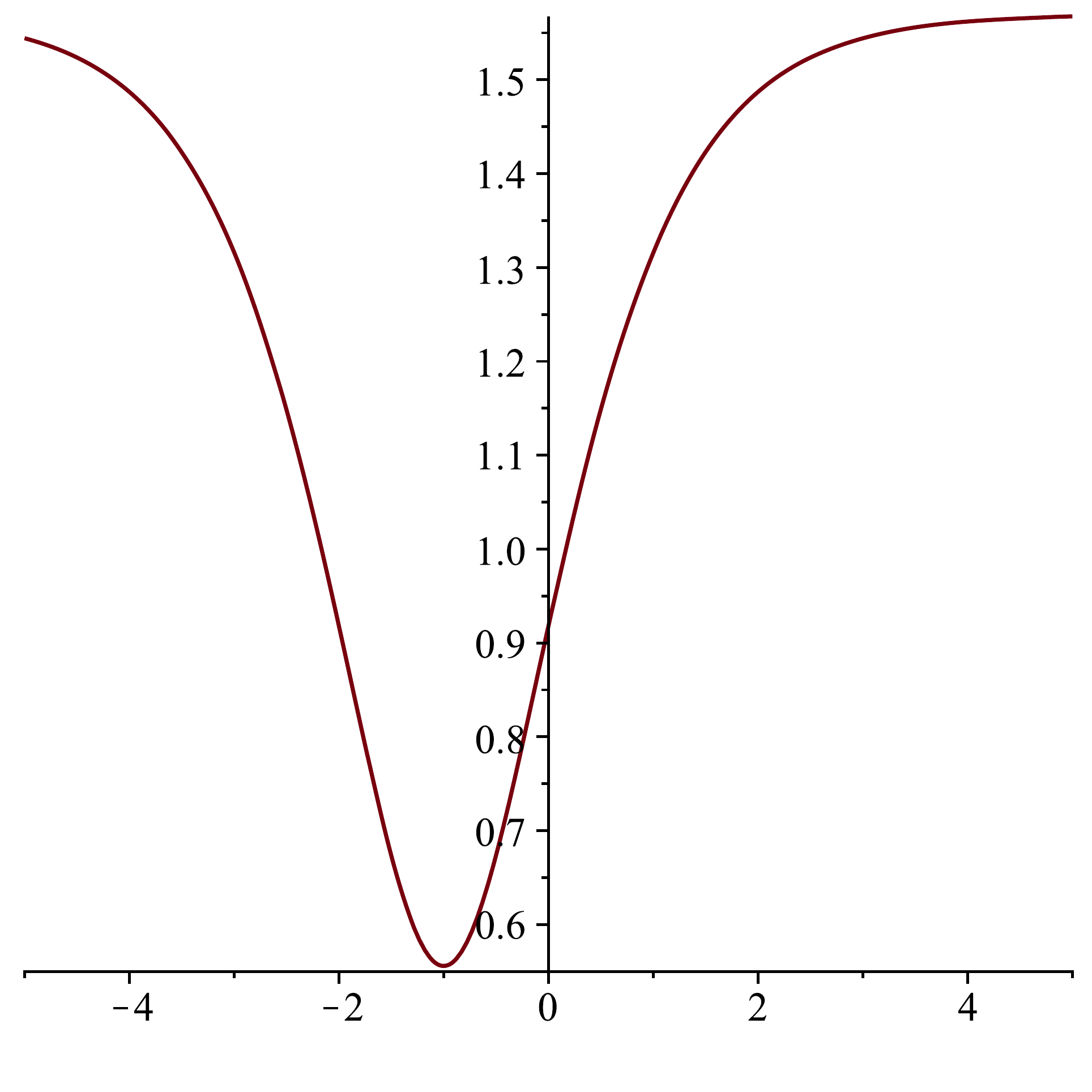}
 		\subcaption{}
 	\end{subfigure}	
 	\caption{(a) Potential $\ds V_2(x)$ of equation (\ref{potential}), with the Green function $\frac{\mu}{\sqrt{\lambda}}\tanh(\frac{|x - x_1|~\mu}{\sqrt{2}})$, $\mu=\lambda=1$, convolved with the Gaussian distribution $e^{-x_1^{2}/a^{2}}$, $a=0.01$; and\\ 
 	(b) Potential $\ds V_3(x)$ of eq. (\ref{potential}), with the same Green function as in (a), $\mu=\lambda=1$, convolved with the distribution $\frac{1} {\left[(x_1+a)^2+b^2)^2 \right]}$, $a=b=1$.}
 	\label{fig:potentials2}
 \end{figure}

 The numerically evaluated integrals (\ref{potential}) (with $\mu = \lambda = 1$, and with  $ a = b = 1$), are plotted in Figure \ref{fig:potentials2}.
 \vskip .2 cm 

\begin{figure}[H]
 	\centering
 	\begin{subfigure}{0.48\textwidth}
 		\includegraphics[width=\textwidth]{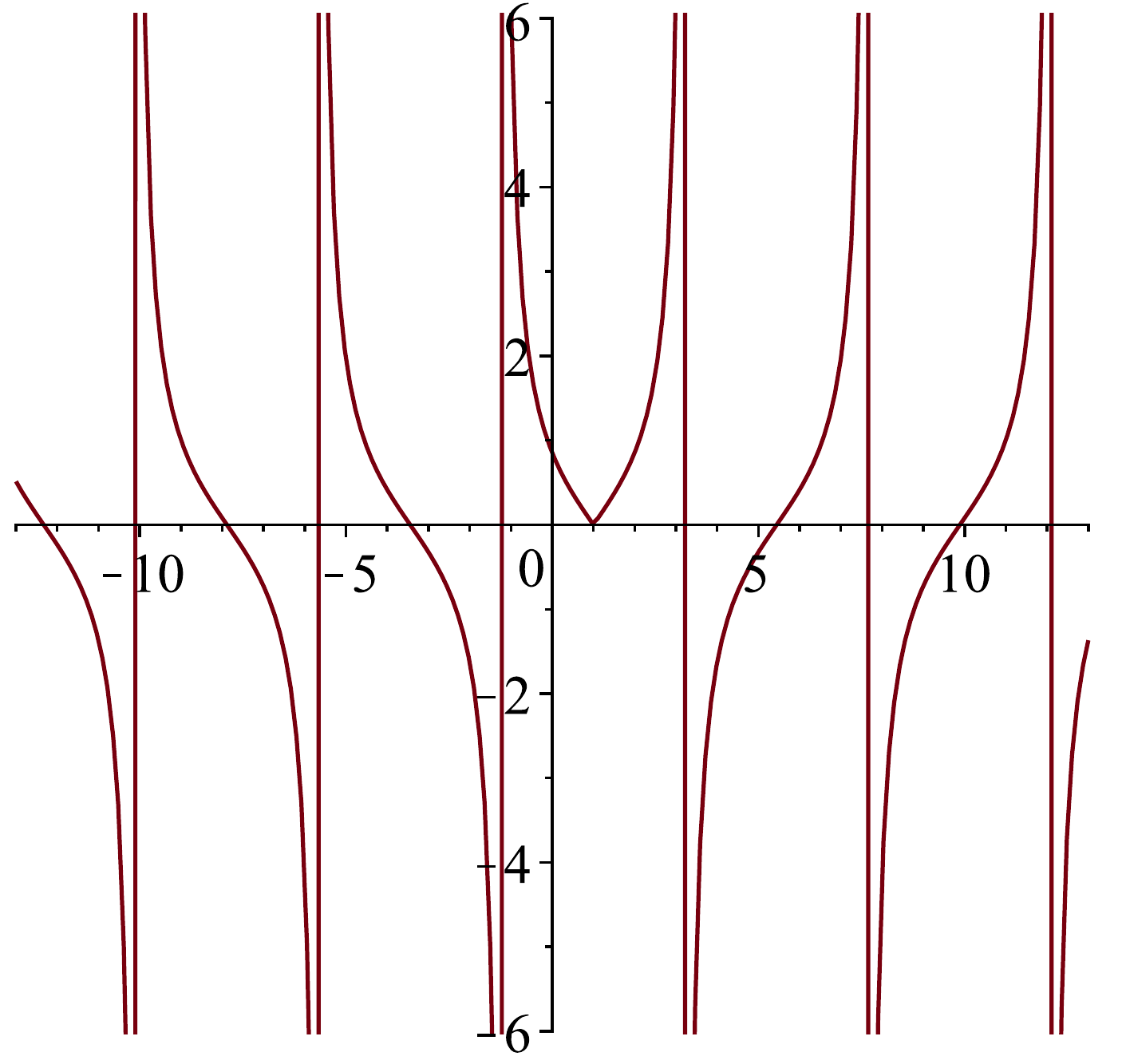}
 		\subcaption{}
 	\end{subfigure}
 	\begin{subfigure}{0.48\textwidth} 
 		\includegraphics[width=\textwidth]{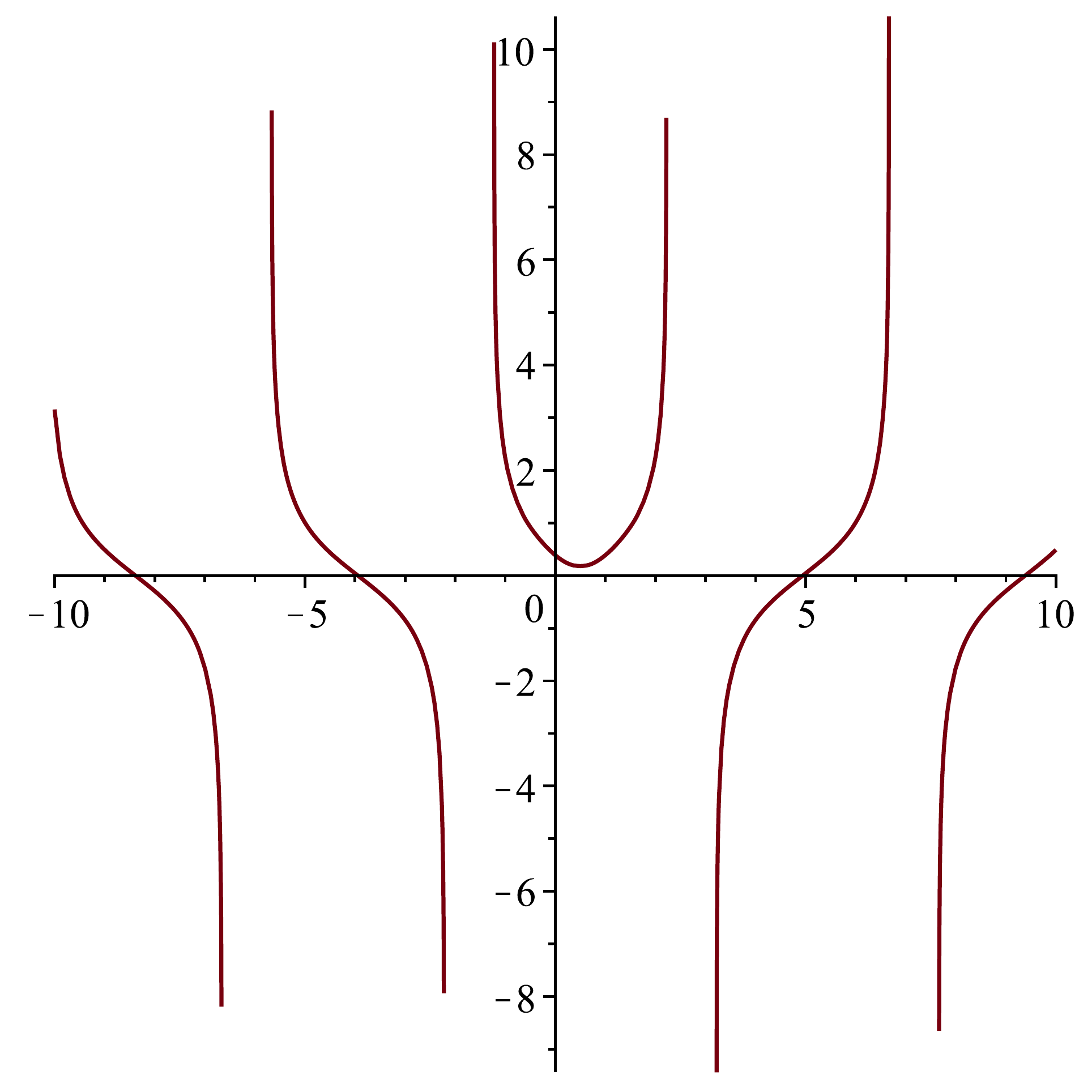}
 		\subcaption{}
 	\end{subfigure}	
 	\caption{(a) Green function $\Psi(x,x_1; m, \lambda) = \frac{m}{ \sqrt\lambda} \tan \left(\frac{|x-x_1| m}{\sqrt2} \right)$, for $x_1=1$, $ m =\lambda=1$  and\ (b) Potential $V(x)$ of the Green function  $\Psi(x,x_1; m=1, \lambda=1)$ convolved with the ``step function'' distribution $R_5(x_1)=1$ $\forall$ $0 \le x_1 \le 1$; $\mu=\lambda=1$.}
 	\label{fig:Greenfunctions2}
 \end{figure}

 The Green function given in eq. (\ref{Greenfunctionpsi}) is quite different from that of eq. (\ref{20}), in that it has an infinity of singular points (vertical asymptotes), as indicated in figure \ref{fig:Greenfunctions2}a 
\vskip .2 cm \ni

Nevertheless it is possible to work out potentials in specified segments between adjacent vertical asymptotes. For example, if we consider a step-function distribution  
 $R_5(x_1) = 1$  for $0 \le x_1 \le 1$ and otherwise $R_5 = 0$; 
 then the potential can be worked out analytically. Its plot is shown in figure \ref{fig:Greenfunctions2}b.
 
\begin{figure}[H]
	\centering
	\begin{subfigure}{0.48\textwidth}
		\includegraphics[width=\textwidth]{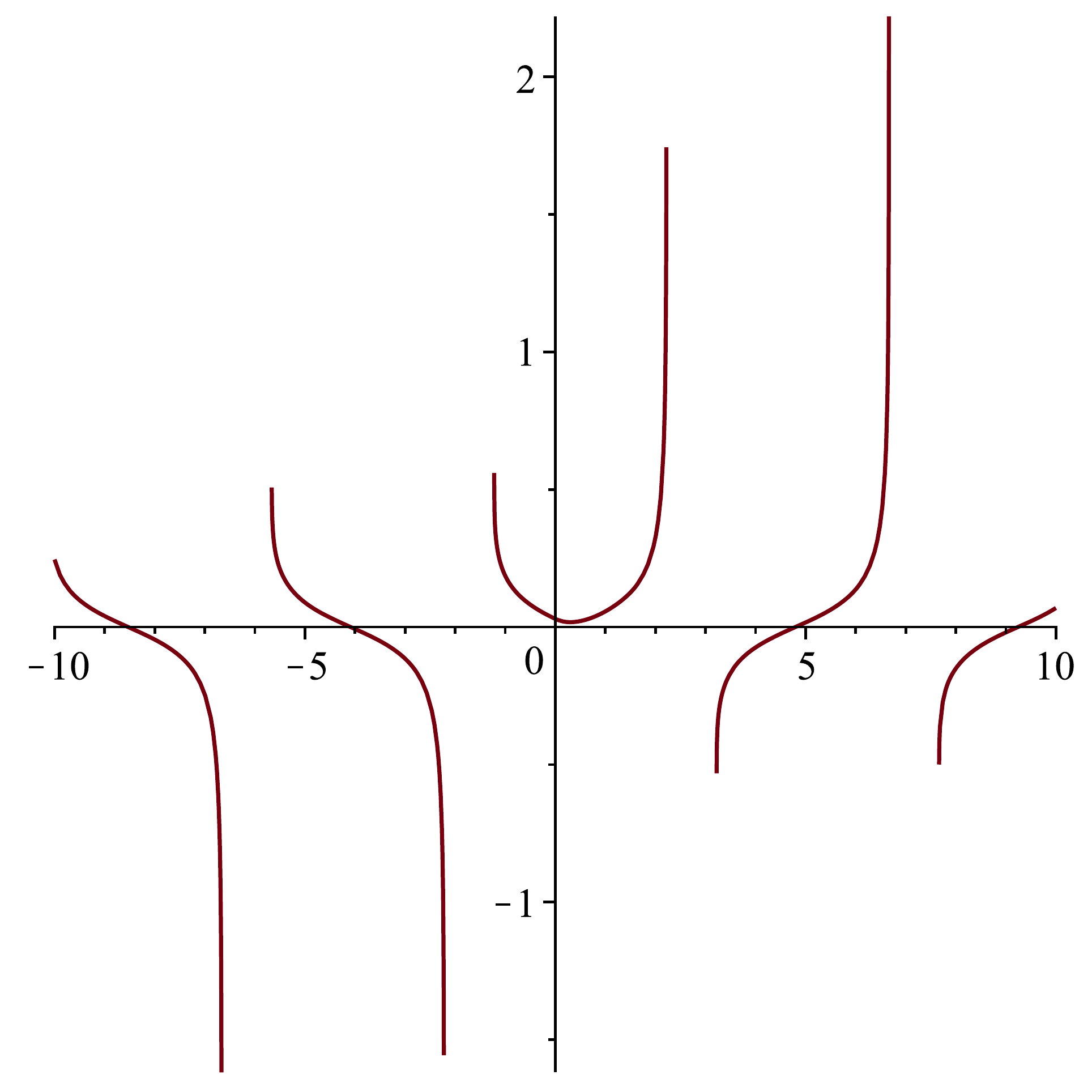}
		\subcaption{}
	\end{subfigure}
	\begin{subfigure}{0.48\textwidth} 
		\includegraphics[width=\textwidth]{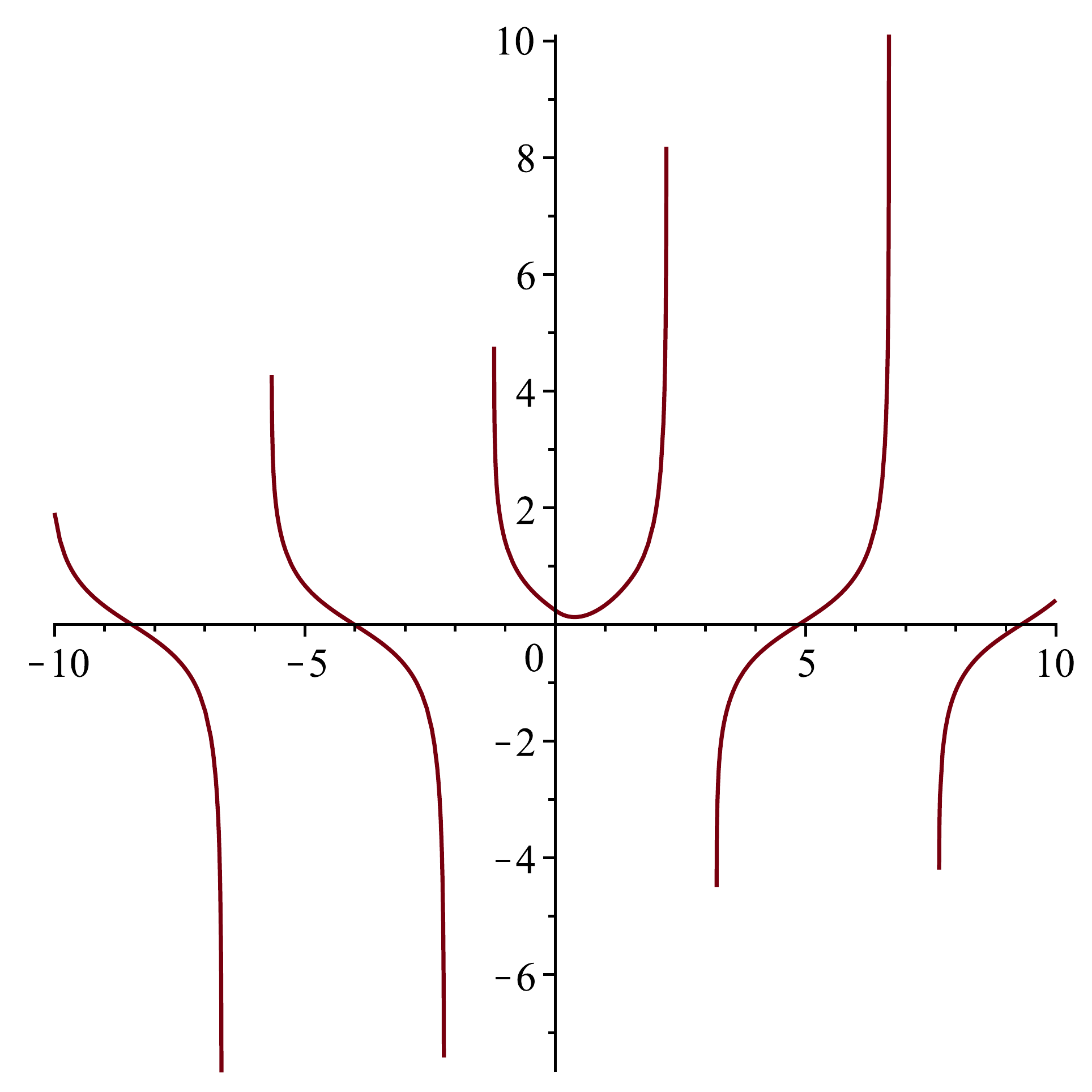}
		\subcaption{}
	\end{subfigure}	
	\caption{(a) Potential from the Green function $\frac{m}{ \sqrt{\lambda}} \tan \left(\frac{|x-x_1| m}{\sqrt 2 } \right)$,
	$m = \lambda=1$,  convolved with the bell shaped distribution $R_3$ (eq. \ref{bellsourcedist}),  with $a=b=1$, and\\ 
	(b) Potential from the  Green function $\frac{m}{\sqrt{\lambda}} \tan \left(\frac{|x-x_{1}| m}{\sqrt{2}}\right)$, $m = \lambda = 1$,  convolved with the Gaussian distribution $R_2$ (eq. \ref{gaussiansourcedist})  with  $a=1$.}
	\label{fig:potentialGreenfunctions2}
\end{figure}

 \vskip .2 cm 
\vskip .2 cm 
For the Gaussian or bell-shaped distributions given in eqs. (\ref{gaussiansourcedist}) and (\ref{bellsourcedist}) the potentials are as in eq. (\ref{23}) but with the Green function as given  in eq. (\ref{Greenfunctionpsi}). In these cases the integrations must be done numerically in specified segments. We plot two examples of such numerical determination of potentials in Figure \ref{fig:potentialGreenfunctions2}: 
 \vskip .2 cm

\vskip .2 cm 
 \section{Concluding remarks}
We have shown that non-linear second-order equations in one dimension for which exact, analytic solutions can be obtained (cf. eqs. (\ref{1}) and (\ref{7})) can be generalised  to 
 inhomogeneous equations with delta function ``point" sources for which Green's functions can be determined (see eqs. (\ref{20}) and (\ref{Greenfunctionpsi})). If the sources are not point sources (``point charges") but continuous source distributions, then integral summations, that is ``potentials", can be defined (cf. eq. (\ref{23})). These can be evaluated analytically in some cases, and otherwise by numerical quadrature, as discussed and illustrated in section 4.
 \vskip .4cm  

\par
\section{Appendix I}
The potential $\ds V_1(x) = \int_{-\infty}^{\infty} dx_1 \, R_1(x_1) \, G(x, x_1,\mu=1,\lambda=1)$, for the step-function
$R_1(x_1) = 1 \;\;\forall\;\;\{-4 \le x_1 \le 8\}, \;\;R_1 = 0\;\;\forall \{x_1>8,\;x_1<-4\}$ 
evaluated by Mathematica\textsuperscript{\textregistered} is shown in (\ref{MMaOutputvxstepfuncint}).
\begin{equation} \label{MMaOutputvxstepfuncint}
V_1(x) = 
\begin{cases} 
{\ds \sqrt{2} \log \left(\frac{e^{\sqrt{2} (x-8)}+1}{e^{\sqrt{2} (x+4)}+1}\right)+12} & \;\;\; x \le -4 \\
{\ds \sqrt{2} \log \left(\frac{e^{\sqrt{2} (x+4)}+1}{e^{\sqrt{2} (x-8)}+1}\right)-12} & \;\;\; x \ge 8 \\
{\ds \sqrt{2} \log \left[\frac{1}{4} \left(e^{\sqrt{2} (x-8)}+1\right) \left(e^{\sqrt{2} (x+4)}+1\right)\right]-2 (x-2)} & \;\;\; -4 \le x \le 8
\end{cases}
\end{equation}
$\ds V_e(x) = \int_{-\infty}^{\infty} dx_1 \, R_e(x_1) \; G(x, x_1,1,1)\;$ for the exponential distribution $\ds R_e(x) = e^{-|x|}$ is evaluated in Mathematica in terms of Hypergeometric, and Polygamma functions, shown in eq. \ref{convolexponential}.
\begin{multline} \label{convolexponential}
    V_{e}(x)= 
\begin{cases}
-\frac{e^{-\sqrt{2} x}}{\sqrt{2}+2} \bigg[ 
2 \sqrt{2} e^{\sqrt{2} x}-2 \sqrt{2} e^{\sqrt{2} x+x}+4 e^{\sqrt{2} x}-4 e^{\sqrt{2} x+x}\\
-8 e^{\sqrt{2} x} \, _2F_1\left(1,-\frac{1}{\sqrt{2}};1-\frac{1}{\sqrt{2}};-e^{\sqrt{2} x}\right)-\\
4 \sqrt{2} e^{\sqrt{2} x} \, _2F_1\left(1,-\frac{1}{\sqrt{2}};1-\frac{1}{\sqrt{2}};-e^{\sqrt{2} x}\right)\\
-2 e^{\sqrt{2} x} \, _2F_1\left(1,\frac{1}{\sqrt{2}};1+\frac{1}{\sqrt{2}};-e^{-\sqrt{2} x}\right)-\\
\sqrt{2} e^{\sqrt{2} x} \, _2F_1\left(1,\frac{1}{\sqrt{2}};1+\frac{1}{\sqrt{2}};-e^{-\sqrt{2} x}\right)\\
-2 e^{\sqrt{2} x} \, _2F_1\left(1,\frac{1}{\sqrt{2}};1+\frac{1}{\sqrt{2}};-e^{\sqrt{2} x}\right)-\\
\sqrt{2} e^{\sqrt{2} x} \, _2F_1\left(1,\frac{1}{\sqrt{2}};1+\frac{1}{\sqrt{2}};-e^{\sqrt{2} x}\right)+\\
\sqrt{2} \, _2F_1\left(1,1+\frac{1}{\sqrt{2}};2+\frac{1}{\sqrt{2}};-e^{-\sqrt{2} x}\right)+\\
\sqrt{2} e^{2 \sqrt{2} x} \, _2F_1\left(1,1+\frac{1}{\sqrt{2}};2+\frac{1}{\sqrt{2}};-e^{\sqrt{2} x}\right)\\
+ 2 e^{\sqrt{2} x+x} \psi ^{(0)}\left(-\frac{1}{2 \sqrt{2}}\right)+
-2 \sqrt{2} e^{\sqrt{2} x+x} \psi ^{(0)}\left(\frac{1}{4} \left(2-\sqrt{2}\right)\right)-\\
2 e^{\sqrt{2} x+x} \psi ^{(0)}\left(\frac{1}{4} \left(2-\sqrt{2}\right)\right)+2 \sqrt{2} e^{\sqrt{2} x+x} \psi ^{(0)}\left(-\frac{1}{2 \sqrt{2}}\right)
\bigg] & \;\;\;x < 0\\
 \frac{e^{-\sqrt{2} x-x}}{\sqrt{2}+2} \bigg[    
-2 e^{\sqrt{2} x}-\sqrt{2} e^{\sqrt{2} x}+2 e^{\sqrt{2} x+x} \, _2F_1\left(1,\frac{1}{\sqrt{2}};1+\frac{1}{\sqrt{2}};-e^{-\sqrt{2} x}\right)\\
+ \sqrt{2} e^{\sqrt{2} x+x} \, _2F_1\left(1,\frac{1}{\sqrt{2}};1+\frac{1}{\sqrt{2}};-e^{-\sqrt{2} x}\right)-\\
\sqrt{2} e^x \, _2F_1\left(1,1+\frac{1}{\sqrt{2}};2+\frac{1}{\sqrt{2}};-e^{-\sqrt{2} x}\right)-\\
e^{\sqrt{2} x} \psi ^{(0)}\left(\frac{1}{2 \sqrt{2}}\right)-\sqrt{2} e^{\sqrt{2} x} \psi ^{(0)}\left(\frac{1}{2 \sqrt{2}}\right)\\
+e^{\sqrt{2} x} \psi ^{(0)}\left(\frac{1}{4} \left(\sqrt{2}+2\right)\right)+\sqrt{2} e^{\sqrt{2} x} \psi ^{(0)}\left(\frac{1}{4} \left(\sqrt{2}+2\right)\right)
\bigg], & \;\;\; x = 0\\
\frac{e^{-\sqrt{2} x-x}}{\sqrt{2}+2} \bigg[ 
-4 e^{\sqrt{2} x}-2 \sqrt{2} e^{\sqrt{2} x}+2 e^{\sqrt{2} x+x} \, _2F_1\left(1,\frac{1}{\sqrt{2}};1+\frac{1}{\sqrt{2}};-e^{-\sqrt{2} x}\right)\\
+ \sqrt{2} e^{\sqrt{2} x+x} \, _2F_1\left(1,\frac{1}{\sqrt{2}};1+\frac{1}{\sqrt{2}};-e^{-\sqrt{2} x}\right)-\\
2 e^{\sqrt{2} x+x} \, _2F_1\left(1,\frac{1}{\sqrt{2}};1+\frac{1}{\sqrt{2}};-e^{\sqrt{2} x}\right)-\\
\sqrt{2} e^{\sqrt{2} x+x} \, _2F_1\left(1,\frac{1}{\sqrt{2}};1+\frac{1}{\sqrt{2}};-e^{\sqrt{2} x}\right)-\\
\sqrt{2} e^x \, _2F_1\left(1,1+\frac{1}{\sqrt{2}};2+\frac{1}{\sqrt{2}};-e^{-\sqrt{2} x}\right)+\\
\sqrt{2} e^{2 \sqrt{2} x+x} \, _2F_1\left(1,1+\frac{1}{\sqrt{2}};2+\frac{1}{\sqrt{2}};-e^{\sqrt{2} x}\right)-\\
2 e^{\sqrt{2} x} \psi ^{(0)}\left(\frac{1}{2 \sqrt{2}}\right)-2 \sqrt{2} e^{\sqrt{2} x} \psi ^{(0)}\left(\frac{1}{2 \sqrt{2}}\right) + \\
2 e^{\sqrt{2} x} \psi ^{(0)}\left(\frac{1}{4} \left(\sqrt{2}+2\right)\right)+2 \sqrt{2} e^{\sqrt{2} x} \psi ^{(0)}\left(\frac{1}{4} \left(\sqrt{2}+2\right)\right)
\bigg], & \;\;\; x >0
\end{cases}
\end{multline}

 %

\end{document}